\documentclass[12 pt]{amsart}
\usepackage{amssymb,bezier}
\newtheorem{theorem}{Theorem}
\newtheorem{definition}[theorem]{Definition}
\newtheorem{proposition}[theorem]{Proposition}
\newtheorem{corollary}[theorem]{Corollary}
\newtheorem{exercise}[theorem]{Exercise}
\newtheorem{example}[theorem]{Example}
\newtheorem{remark}[theorem]{Remark}

\def\Hom{{\it Hom}}                    \def\ot{\otimes}
\def\squeezedldots{{\mbox{$. \hskip -1pt . \hskip -1pt .$}}}
\def\as{{\it as}}                      \def\id{{\mbox{1 \hskip -8pt 1}}}
\def\ext{\mbox{\large$\land$}}         \def\bbL{{\mathbb{L}}}
\def\coland{\mbox{${}^c\!\ext$}}       \def\angles#1{{\langle #1 \rangle}}
\def\susp{\uparrow\!}                  \def\sb#1{\langle #1 \rangle}
\def\desusp{\downarrow\!}
\def\rada#1#2#3{{#1}_{#2},\ldots,{#1}_{#3}}
\def\Rada#1#2#3{#1_{#2},\dots,#1_{#3}}
\def\sgn{\mbox{\rm sgn}}
\def\otexp#1#2{{#1}^{\otimes #2}}

\def\Alpha#1{{
{
\thicklines
\unitlength=#1pt
\begin{picture}(100.00,110.00)(0.00,0.00)
\put(50.00,50.00){\makebox(0.00,0.00){\scriptsize $\alpha$}}
\put(65.00,15.00){\makebox(0.00,0.00){$\cdots$}}
\put(50.00,50.00){\line(0,1){0.00}}
\put(50.00,110.00){\line(0,-1){30.00}}
\put(50.00,95.00){\makebox(0,0){$\bullet$}}
\put(90.00,30.00){\line(0,-1){30.00}}
\put(30.00,30.00){\line(0,-1){30.00}}
\put(10.00,30.00){\line(0,-1){30.00}}
\put(100.00,30.00){\bezier{50}(0,0)(-25,25)(-50,50)}
\put(0.00,30.00){\line(1,0){100.00}}
\put(50,80){\bezier{50}(0,0)(-25,-25)(-50,-50)}
\end{picture}}
}}

\def\Alphanobul#1{{
{
\thicklines
\unitlength=#1pt
\begin{picture}(100.00,110.00)(0.00,0.00)
\put(50.00,50.00){\makebox(0.00,0.00){\scriptsize $\alpha$}}
\put(65.00,15.00){\makebox(0.00,0.00){\multiput(-20,0)(10,0){3}{$\cdot$}}}
\put(50.00,50.00){\line(0,1){0.00}}
\put(50.00,110.00){\line(0,-1){30.00}}
\put(90.00,30.00){\line(0,-1){30.00}}
\put(30.00,30.00){\line(0,-1){30.00}}
\put(10.00,30.00){\line(0,-1){30.00}}
\put(100.00,30.00){\bezier{50}(0,0)(-25,25)(-50,50)}
\put(0.00,30.00){\line(1,0){100.00}}
\put(50,80){\bezier{50}(0,0)(-25,-25)(-50,-50)}
\end{picture}}
}}

\def\Beta#1{{
{
\thicklines
\unitlength=#1pt
\begin{picture}(100.00,110.00)(0.00,0.00)
\put(50.00,50.00){\makebox(0.00,0.00){\scriptsize $\beta$}}
\put(65.00,15.00){\makebox(0.00,0.00){\multiput(-20,0)(10,0){3}{$\cdot$}}}
\put(50.00,50.00){\line(0,1){0.00}}
\put(50.00,110.00){\line(0,-1){30.00}}
\put(50.00,95.00){\makebox(0,0){$\bullet$}}
\put(90.00,30.00){\line(0,-1){30.00}}
\put(30.00,30.00){\line(0,-1){30.00}}
\put(10.00,30.00){\line(0,-1){30.00}}
\put(100.00,30.00){\bezier{50}(0,0)(-25,25)(-50,50)}
\put(0.00,30.00){\line(1,0){100.00}}
\put(50,80){\bezier{50}(0,0)(-25,-25)(-50,-50)}
\end{picture}}
}}

\def\Gamma#1{{
{
\thicklines
\unitlength=#1pt
\begin{picture}(100.00,110.00)(0.00,0.00)
\put(50.00,50.00){\makebox(0.00,0.00){\scriptsize $\gamma$}}
\put(65.00,15.00){\makebox(0.00,0.00){\multiput(-20,0)(10,0){3}{$\cdot$}}}
\put(50.00,50.00){\line(0,1){0.00}}
\put(50.00,110.00){\line(0,-1){30.00}}
\put(50.00,95.00){\makebox(0,0){$\bullet$}}
\put(90.00,30.00){\line(0,-1){30.00}}
\put(30.00,30.00){\line(0,-1){30.00}}
\put(10.00,30.00){\line(0,-1){30.00}}
\put(100.00,30.00){\bezier{50}(0,0)(-25,25)(-50,50)}
\put(0.00,30.00){\line(1,0){100.00}}
\put(50,80){\bezier{50}(0,0)(-25,-25)(-50,-50)}
\end{picture}}
}}

\def\Gammanobul#1{{
{
\thicklines
\unitlength=#1pt
\begin{picture}(100.00,110.00)(0.00,0.00)
\put(50.00,50.00){\makebox(0.00,0.00){\scriptsize $\gamma$}}
\put(65.00,15.00){\makebox(0.00,0.00){\multiput(-20,0)(10,0){3}{$\cdot$}}}
\put(50.00,50.00){\line(0,1){0.00}}
\put(50.00,110.00){\line(0,-1){30.00}}
\put(90.00,30.00){\line(0,-1){30.00}}
\put(30.00,30.00){\line(0,-1){30.00}}
\put(10.00,30.00){\line(0,-1){30.00}}
\put(100.00,30.00){\bezier{50}(0,0)(-25,25)(-50,50)}
\put(0.00,30.00){\line(1,0){100.00}}
\put(50,80){\bezier{50}(0,0)(-25,-25)(-50,-50)}
\end{picture}}
}}

\def\Betanobul#1{{
{
\thicklines
\unitlength=#1pt
\begin{picture}(100.00,110.00)(0.00,0.00)
\put(50.00,50.00){\makebox(0.00,0.00){\scriptsize $\beta$}}
\put(65.00,15.00){\makebox(0.00,0.00){\multiput(-20,0)(10,0){3}{$\cdot$}}}
\put(50.00,50.00){\line(0,1){0.00}}
\put(50.00,110.00){\line(0,-1){30.00}}
\put(90.00,30.00){\line(0,-1){30.00}}
\put(30.00,30.00){\line(0,-1){30.00}}
\put(10.00,30.00){\line(0,-1){30.00}}
\put(100.00,30.00){\bezier{50}(0,0)(-25,25)(-50,50)}
\put(0.00,30.00){\line(1,0){100.00}}
\put(50,80){\bezier{50}(0,0)(-25,-25)(-50,-50)}
\end{picture}}
}}

\def\NABLA#1{{
{
\thicklines
\unitlength=#1pt
\begin{picture}(100.00,110.00)(0.00,0.00)
\put(50.00,50.00){\makebox(0.00,0.00){\scriptsize $\nabla$}}
\put(65.00,15.00){\makebox(0.00,0.00){\multiput(-20,0)(10,0){3}{$\cdot$}}}
\put(50.00,50.00){\line(0,1){0.00}}
\put(50.00,110.00){\line(0,-1){30.00}}
\put(90.00,30.00){\line(0,-1){30.00}}
\put(30.00,30.00){\line(0,-1){30.00}}
\put(10.00,30.00){\line(0,-1){30.00}}
\put(10.00,15.00){\makebox(0,0){$\bullet$}}
\put(100.00,30.00){\line(-1,1){50.00}}
\put(0.00,30.00){\line(1,0){100.00}}
\put(50.00,80.00){\line(-1,-1){50.00}}
\end{picture}}
}}

\def\NABLAmod#1{{
{
\thicklines
\unitlength=#1pt
\begin{picture}(100.00,110.00)(0.00,0.00)
\put(50.00,50.00){\makebox(0.00,0.00){\scriptsize $\nabla$}}
\put(65.00,15.00){\makebox(0.00,0.00){$\cdots$}}
\put(50.00,50.00){\line(0,1){0.00}}
\put(50.00,110.00){\line(0,-1){30.00}}
\put(90.00,30.00){\line(0,-1){30.00}}
\put(35.00,30.00){\line(0,-1){30.00}}
\put(35.00,15.00){\makebox(0,0){$\bullet$}}
\put(5.00,30.00){\line(0,-1){30.00}}
\put(100.00,30.00){\line(-1,1){50.00}}
\put(0.00,30.00){\line(1,0){100.00}}
\put(50.00,80.00){\line(-1,-1){50.00}}
\end{picture}}
}}

\def\NABLAnobul#1{{
{
\thicklines
\unitlength=#1pt
\begin{picture}(100.00,110.00)(0.00,0.00)
\put(52.00,47.00){\makebox(0.00,0.00){\scriptsize $\nabla$}}
\put(65.00,15.00){\makebox(0.00,0.00){\multiput(-20,0)(10,0){3}{$\cdot$}}}
\put(50.00,50.00){\line(0,1){0.00}}
\put(50.00,110.00){\line(0,-1){30.00}}
\put(90.00,30.00){\line(0,-1){30.00}}
\put(30.00,30.00){\line(0,-1){30.00}}
\put(10.00,30.00){\line(0,-1){30.00}}
\put(100.00,30.00){\bezier{50}(0,0)(-25,25)(-50,50)}
\put(0.00,30.00){\line(1,0){100.00}}
\put(50,80){\bezier{50}(0,0)(-25,-25)(-50,-50)}
\end{picture}}
}}

\def\UPSILON#1{{
{
\thicklines
\unitlength=#1pt
\begin{picture}(100.00,110.00)(0.00,0.00)
\put(50.00,50.00){\makebox(0.00,0.00){\scriptsize $\Upsilon$}}
\put(65.00,15.00){\makebox(0.00,0.00){$\cdots$}}
\put(50.00,50.00){\line(0,1){0.00}}
\put(50.00,110.00){\line(0,-1){30.00}}
\put(50.00,95.00){\makebox(0,0){$\bullet$}}
\put(90.00,30.00){\line(0,-1){30.00}}
\put(5,0){
\put(30.00,30.00){\line(0,-1){30.00}}
\put(30.00,15.00){\makebox(0,0){$\bullet$}}
}
\put(-5,0){
\put(10.00,30.00){\line(0,-1){30.00}}
\put(10.00,15.00){\makebox(0,0){$\bullet$}}
}
\put(100.00,30.00){\line(-1,1){50.00}}
\put(0.00,30.00){\line(1,0){100.00}}
\put(50.00,80.00){\line(-1,-1){50.00}}
\end{picture}}
}}

\begin{document}
\baselineskip 17pt
\bibliographystyle{alpha}

\title{Symmetric Brace Algebras}
\author{Tom Lada}
\address{Department of Mathematics, North Carolina State University,
Raleigh NC 27695, USA}
\email{lada@math.ncsu.edu}
\thanks{The research of the first author was supported in part by NSF
grant INT-0203119}
\author{Martin Markl}
\address{Mathematics Institute of the Academy, \v Zitn\'a 25, 115 67 Praha 1,
The Czech Republic}
\email{markl@math.cas.cz}
\thanks{The research of the second author was supported by 
grant M{\v S}MT ME 603  and by
the Academy of Sciences of the Czech Republic, 
Institutional Research Plan No.~AV0Z10190503.}

\begin{abstract}
We develop a symmetric analog of brace algebras and discuss the relation
of such algebras to $L_{\infty}$-algebras. 
We give an alternate proof  that the category of symmetric brace algebras is isomorphic to the
category of pre-Lie algebras.  As an application, symmetric braces are used to describe
transfers of strongly homotopy structures.
We then explain how these symmetric 
brace algebras
may be used to examine the $L_{\infty}$-algebras that result from a particular
gauge theory for massless particles of high spin.
\end{abstract}

\maketitle

\parskip 4pt

\section*{Introduction}

The interplay between homological algebra and mathematical physics
provides a vast arena for research in both subjects.  In algebra, the
construction of multilinear operations, called braces,  on the
Hochschild complex of an associative algebra leads to the definition of
a brace algebra  on a graded vector space~\cite{Kad} 
and to subsequent applications in topological
field theory \cite{KVZ}, as well as in a particular solution to the Deligne conjecture \cite{G,GV}.   In this note, we
develop the idea of a {\em symmetric brace algebra\/} in which the brace
operations will possess the property of graded symmetry.  Multilinear operations
on the space of anti-symmetric maps on a graded vector space, e.g. the
Chevalley-Eilenberg complex of a Lie algebra, provides a motivating example 
for this concept.  Just as $A_{\infty}$- algebraic relations may be neatly encoded in the brace algebra context, $L_{\infty}$-algebra data is given by symmetric brace algebra relations.

Both brace and symmetric brace algebra structures yield pre-Lie algebra structures on the underlying vector space.  However, a major difference is that a pre-Lie algebra may be used to define a symmetric brace algebra structure \cite{guin-oudom}.  We provide a conceptual proof of this fact and along the way exhibit a construction of a model of a free symmetric algebra.

As a first application of symmetric braces, suppose that we are given chain complexes $(V,\partial_V)$, $(W,\partial_W)$ and chain maps
$f : (V,\partial_V) \to (W,\partial_W)$, $g : (W,\partial_W) \to (V,\partial_V)$ such
that the composition $gf$ is chain homotopic to the identity $\id_V:
V\to V$, via a chain-homotopy $h$.  If $V$ possesses an $A_\infty$- structure or an $L_\infty$-structure,  we show that braces and symmetric braces respectively, may be used to transfer these homotopy structures to the space $W$.

As another application of symmetric braces, we recall some work of Berends, Burgers, and van~ Dam in mathematical physics that generated algebraic data describing interactions between massless particles of high spin \cite{BBvD}.  
In \cite{FLS} it was shown that these algebraic relations may be reformulated
as giving the graded vector space of gauge parameters together with the
fields the structure of an $L_{\infty}$-algebra.  In this note, we show
that this $L_{\infty}$-structure may be conveniently analyzed using
the language of symmetric brace algebras. This places the problem into a more
general algebraic context which may perhaps be applicable to other
similar problems as well.

These applications demonstrate that even though the category of symmetric brace algebras is isomorphic to  the category of pre-Lie algebras, symmetric braces may provide a more manageable setting for calculations than pre-Lie operations.

Each (nonsymmetric) brace algebra defines a symmetric one by an
obvious symmetrization, but not every symmetric algebra is the
symmetrization of a nonsymmetric one. In this sense, nonsymmetric 
braces are more
special than symmetric ones. While symmetric braces
can be expected to exist on the operadic cohomology complex
of any algebra over a Koszul quadratic operad, nonsymmetric braces exist
only for algebras over non-$\Sigma$ operads,
compare also the remarks in~\cite[II.3.9]{markl-shnider-stasheff:book}.

We recall, in Section~\ref{sec2}, the definition of brace algebras  
and develop our
concept of symmetric brace algebras.  We discuss the relation of such
algebras to $A_{\infty}$ and to $L_{\infty}$-algebras.  Several general
properties of symmetric brace algebras are also given here.
Section~\ref{sec3} contains a proof that the category of symmetric brace algebras is isomorphic to the category of pre-Lie algebras.  The transfer of strongly homotopy structures appears in Section ~\ref{sec4}.
In Section~\ref{sec5}, we review the algebraic background 
of the physics problem
as formulated in~\cite{BBvD,FLS}.  We place this problem into the context of
symmetric brace algebras and identify several relations in this algebra
that yield the same $L_{\infty}$-structure that is found in \cite{FLS}.

\section{Symmetric Brace Algebras}
\label{sec2}

We begin by introducing the following necessary technical
notions~\cite[page~2148]{LM}. For graded indeterminates 
$\rada x1n$ and a permutation $\sigma\in
\Sigma_n$ define the {\em Koszul sign\/}
$\epsilon =\epsilon(\sigma;\rada x1n)$ by
\[
x_1\land\dots\land x_n = \epsilon(\sigma;x_1,\dots,x_n)
\cdot x_{\sigma(1)}\land \dots \land x_{\sigma(n)},
\]
which has to be satisfied in the free graded commutative algebra
$\ext(\rada x1n)$. We will need also the {\em antisymmetric Koszul sign\/}
\[
\chi = \chi(\sigma;\rada x1n)
:= \sgn(\sigma)\cdot \epsilon(\sigma;\rada x1n).
\]
Let us recall the definition of a (non-symmetric) brace
algebra as given in~\cite{GV}.

\begin{definition}  A brace algebra is a graded vector space $U$ together
with a collection of degree~$0$ multilinear braces
$x,x_1,\ldots,x_n \longmapsto  x\{x_1,\ldots,x_n\}$  that
satisfy the identities
\[
x\{\thickspace\}=x
\]
and
\begin{eqnarray*}
\lefteqn{
x\{x_1,\squeezedldots,x_m\}\{y_1,\squeezedldots,y_n\}=}
\\
&&
\sum \epsilon \cdot x\{y_1,\squeezedldots,y_{i_1},x_1\{y_{i_1+1},\squeezedldots,y_{j_1}\},
y_{j_1+1},\squeezedldots,y_{i_m},
x_m\{y_{i_m+1},\squeezedldots,y_{j_m}\},y_{j_m+1},\squeezedldots,y_n\},
\end{eqnarray*}
where the sum is taken over all sequences 
$0\leq i_1\leq j_1 \leq \cdots\leq i_m \leq j_m  \leq n$ 
and where $\epsilon$ is the Koszul sign of the permutation
\[
(x_1,\squeezedldots,x_m, y_1,\squeezedldots,y_n) \mapsto
(y_1,\squeezedldots,y_{i_1},x_1,y_{i_1+1},\squeezedldots,y_{j_1},y_{j_1+1},
\squeezedldots,y_{i_m},x_m,y_{i_m+1},\squeezedldots,y_{j_m},y_{j_m+1},
\squeezedldots,y_n)
\]
of elements of $U$.
\end{definition}

In~\cite{GV}, two gradings of the underlying vector space, $\deg(x)$ and $|x|$,
related by $|x| = \deg(x)-1$,  were used. 
The brace $x\{x_1,\ldots,x_n\}$ there was of
degree $-n$ with respect to the $\deg(-)$-grading and of degree $0$ with
respect to the $|\hskip -1mm- \hskip -1mm|$-grading. 
When we refer to~\cite{GV},
we always consider the underlying vector space graded with the
$|\hskip -1mm - \hskip -1mm|$-grading. Then
all braces will be degree zero maps, as assumed in the above definition.
We now consider a symmetric version of the brace algebra.

\begin{definition}
\label{sb}
A symmetric brace algebra is a graded vector space $B$ together with a
collection of degree $0$ multilinear braces
$x\langle x_1,\ldots,x_n\rangle$ that are graded symmetric in
$x_1,\ldots,x_n$ and satisfy the identities
\[
x\langle\thickspace\rangle=x
\]
and
\begin{eqnarray}
\label{pozitri_prohlidka_a_nejsem_zdrav}
\lefteqn{
x\langle x_1,\ldots,x_m\rangle\langle y_1,\ldots,y_n\rangle=}
\\
\nonumber 
&&
\sum\epsilon \cdot x\langle x_1\langle y_{i_1^1},\ldots, y_{i_{t_1}^1}\rangle,
x_2\langle y_{i_1^2},\ldots,y_{i_{t_2}^2}\rangle,\ldots,
x_m\langle y_{i_1^m},\ldots,y_{i_{t_m}^m}\rangle,y_{i_1^{m+1}},\ldots,
y_{i_{t_{m+1}}^{m+1}}\rangle
\end{eqnarray}
where the sum is taken over all unshuffle sequences
$$
i_1^1<\cdots<i_{t_1}^1,\ldots,i_1^{m+1}<\cdots<i_{t_{m+1}}^{m+1}
$$
of $\{1,\ldots,n\}$ and where $\epsilon$ is the Koszul sign of the permutation
\[
(x_1,\squeezedldots,x_m, y_1,\squeezedldots,y_n) \longmapsto
(x_1,  y_{i_1^1},\squeezedldots, y_{i_{t_1}^1},
x_2,y_{i_1^2},\squeezedldots,y_{i_{t_2}^2},\squeezedldots,
x_m, y_{i_1^m},\squeezedldots,y_{i_{t_m}^m},y_{i_1^{m+1}},\squeezedldots,
y_{i_{t_{m+1}}^{m+1}})
\]
of elements of $B$.
\end{definition}

\begin{exercise}
\label{bylo_to_silenstvi}
{\rm
For elements $x,y$ of an arbitrary symmetric brace algebra $B$, put
\[
x \circ y := x \angles y.
\]
Prove that then $(B,\circ)$ is a graded right pre-Lie algebra in the sense
of~\cite[Section~2]{gerstenhaber:AM63}, therefore 
$[x,y] := x \circ y - (-1)^{|x||y|} y \circ x$
defines a graded Lie algebra structure on~$B$. This result should be
compared to a similar result for non-symmetric brace algebras
in~\cite[p.~143]{GV}.

This Lie algebra structure corresponds to the standard anti-commutator
Lie algebra on the space of coderivations of the cocommutative
coassociative coalgebra $\coland(W)$~\cite[page~2150]{LM}
cogenerated by the suspension $W:= \hskip 1mm \susp \hskip -.3mm V$ of $V$. 
Compare also remarks in~\cite[II.3.9]{markl-shnider-stasheff:book}. 
}
\end{exercise}

\begin{remark}
\label{nejak_se_nic_nedari}
{\rm
J.-M. Oudom and D. Guin proved in \cite{ guin-oudom} that higher brackets $x\langle \Rada x1n
\rangle$ of an arbitrary symmetric brace algebra are, 
for $n \geq 2$, determined by the `pre-Lie part' $x \circ y = x\langle
y \rangle$, introduced in Exercise~\ref{bylo_to_silenstvi}. For
instance, axiom~(\ref{pozitri_prohlidka_a_nejsem_zdrav}) 
implies that $x \langle x_1,x_2 \rangle$ can be expressed as 
\[
x \langle x_1,x_2 \rangle =  x \sb {x_1} \sb {x_2} - 
x\sb {x_1 \sb {x_2}} = (x \circ x_1) \circ x_2 - x \circ (x_1
\circ x_2). 
\]
The same axiom applied on $x\sb{\Rada x1{n-1}}\sb{x_n}$ can then be clearly
interpreted as an inductive rule defining $x\sb{ \Rada x1n}$
in terms of $x\langle \Rada x1k \rangle$, with $k < n$. 

They also proved that {\em an arbitrary\/} pre-Lie algebra
determines in this way a symmetric brace algebra, which would mean
that the category of symmetric brace algebras {\em is isomorphic\/} to
the category of pre-Lie algebras. Let us observe that the proof of
this statement is not obvious. First,
axiom~(\ref{pozitri_prohlidka_a_nejsem_zdrav}) interpreted as an
inductive rule is `overdetermined.' For example, $x\sb {x_1,x_2,x_3}$
can be expressed both from~(\ref{pozitri_prohlidka_a_nejsem_zdrav})
applied to $x \sb {x_1,x_2}\sb {x_3}$ and also
from~(\ref{pozitri_prohlidka_a_nejsem_zdrav}) applied to $x \sb {x_1}
\sb {x_2,x_3}$, and it is not obvious that the results are the
same. Second, even if the braces are well-defined, it is not clear
that they satisfy the axioms of brace algebras, including the
graded symmetry. 
}
\end{remark}

\begin{example}
\label{Jsem_zvedav_jestli_se_ozve_pred_Sevillou}
{\rm
Just as there is a  (nonsymmetric) brace algebra structure on the graded
vector space $\bigoplus_{k \geq 1}\Hom(V^{\otimes k},V)$ (see~\cite{GV}), 
the basic example of a
symmetric brace algebra is provided by the space of antisymmetric
(another terminology: alternating)
maps, $\bigoplus_{k \geq 1} \Hom(V^{\otimes k},V)^{\as}$. 
More precisely, let $B(V) = B_*(V)$ be
the graded vector space with components
\[
B_s(V) := \bigoplus_{p-k+1 = s} \Hom(V^{\otimes k},V)^{\as}_p,
\]  
where $\Hom(V^{\otimes k},V)^{\as}_p$ denotes the space of 
$k$-multilinear maps of degree $p$ that are antisymmetric
(or alternating) in the sense that 
\[
f(v_1,\ldots,v_i,v_{i+1},\ldots,v_k) = - (-1)^{\deg(v_i)\deg(v_{i+1})}
f(v_1,\ldots,v_{i+1},v_i,\ldots,v_k),
\]
for any $v_1,\ldots,v_i,v_{i+1},\ldots,v_k \in V$ and $1 \leq i \leq k-1$.

Given $f \in \Hom(V^{\otimes k},V)^{\as}_p$ and $g_i \in \Hom(V^{\otimes
a_i},V)^{\as}_{q_i}$, $1 \leq i \leq n$,
define the symmetric brace
$f\langle g_1,\ldots,g_n\rangle \in \Hom(V^{\otimes
r},V)^{\as}_{p+q_1+\cdots+q_n}$, where $r := a_1 + \cdots +a_n + k -n$, by
\begin{equation}
\label{nikdo}
f\langle g_1,\ldots,g_n\rangle(v_1,\ldots,v_{r}) :=
\sum(-1)^{\delta}\chi \cdot f(g_1\otimes\cdots \otimes  g_n
\otimes \id^{\otimes k-n})(v_{i_1},\ldots,v_{i_{r}})
\end{equation}
with the summation taken over all unshuffles
\[
i_1 < \cdots < i_{a_1}, i_{a_1 +1} < \ldots < i_{a_1+a_2},
\ldots, i_{a_1+\cdots + a_k + 1} < \cdots < i_{r},
\]
of elements of $V$, where $\chi$ is the antisymmetric Koszul 
sign of the permutation
\[
(v_1,\ldots,v_{r}) \longmapsto (v_{i_1},\ldots,v_{i_{r}})
\]
and
\begin{eqnarray*}
\delta &=& (k-1)q_1+(k-2+a_1)q_2+\cdots+(k-n + a_1+\cdots+a_{n-1})q_n
\\
&&+  \textstyle\sum_{1 \leq i < j \leq n}a_ia_j  +
(n-1)a_1+(n-2)a_2+\cdots+a_{n-1}.
\end{eqnarray*}
}
\end{example}

\begin{exercise}
\label{Lowecraft}
{\rm
Just as an $A_{\infty}$-structure on $V$ may be described by the 
brace algebra relation 
$\mu\{\mu\}=0$, with $\mu=\mu_1+\mu_2+\cdots$, 
$\mu_k \in \Hom(V^{\otimes k},V)_{k-2}$, an 
$L_\infty$-algebra structure on $V$ can be described by the
symmetric brace algebra relation $l\langle l\rangle=0$; 
here $l=l_1+l_2+\cdots$ where each $l_k\in
\Hom(V^{\otimes k},V)^{\as}_{k-2} \in B_{-1}(V)$. 
Strictly speaking, elements $\mu$
and $l$ belong to the completions  $\prod_{k \geq 1} 
\Hom(V^{\otimes k},V)_{k-2}$ and $\prod_{k \geq 1} 
\Hom(V^{\otimes k},V)_{k-2}^\as$ of the underlying graded vector
spaces, but it is immediately clear that the
above statements make sense also in this more general setup.
}
\end{exercise}

\begin{exercise}
\label{dnes_se_mi_o_ni_zdalo}
{\rm
As a very particular case of Exercise~\ref{Lowecraft}, each Lie
algebra structure on $V$ determines an element $l = l_2 \in
\Hom(V^{\otimes 2},V)_0^\as \subset B_{-1}(V)$ such that $l\langle l
\rangle = 0$. Prove that then the formulas
\[
\partial f  := l\langle f \rangle - (-1)^{|f|} f \langle l \rangle \mbox
{ and }
\{f,g\} := l \langle f,g \rangle
\]
define on $B(V)$ a differential graded Lie algebra,
with a degree $-1$ bracket $\{-,-\}$ and degree~$-1$ differential
$\partial$. 
Verify also the formula
\[
\{f,g\}  =  \partial f \circ g + (-1)^{|f|}f \circ 
\partial g   -\partial (f \circ g) 
\]
which shows that that the bracket $\{-,-\}$ is
actually {\em cohomologous to zero\/}, with the chain
homotopy given by $f \circ g$.
}
\end{exercise}

It is easy to see that, in the situation of
Exercise~\ref{dnes_se_mi_o_ni_zdalo}, the bigraded complex 
\[
{\it CE}^*(V):= (B_{1-*}(V),\partial)
\]
is the standard Chevalley-Eilenberg complex of the graded Lie algebra
$(V,l)$. Brackets $[-,-]$ and $\{-,-\}$, introduced in
Exercises~\ref{bylo_to_silenstvi} and~\ref{dnes_se_mi_o_ni_zdalo},
induce on ${\it CE}^*(V)$ two Lie brackets, which we denote again by
$[-,-]$ and $\{-,-\}$, of degrees $-1$ and $0$ respectively. The first
bracket $[-,-]$, whose definition does not involve $l$, is the
intrinsic bracket considered, for example,
in~\cite{schlessinger-stasheff:JPAA85}. The second bracket $\{-,-\}$  
should be considered as an
analog of the $\smile$\ -product on Hochschild cochain complex of an
associative algebra,
see~\cite{GV}. The above calculation shows that this bracket is
homologically trivial, therefore one could not expect the
Chevalley-Eilenberg cohomology of Lie algebras 
to have a similar rich structure as the Hochschild cohomology of
associative algebras.

The relationship between brace algebras and symmetric brace 
algebras may be summarized by the following theorems.

\begin{theorem}
\label{A}
Let $-\{-,\cdots,-\}$ be a (non-symmetric) brace algebra structure
on a graded vector space $U$. Then
$$
f\langle g_1,\ldots,g_n\rangle:=
\sum_{\sigma\in \Sigma_n} \epsilon \cdot f\{g_{\sigma(1)},\ldots,
g_{\sigma(n)}\},
$$
where $\epsilon$ denotes the Koszul sign of the permutation
\[
(g_1,\ldots,g_n) \longmapsto (g_{\sigma(1)},\ldots, g_{\sigma(n)}),
\]
gives $U$ the structure of a symmetric brace algebra.
\end{theorem}

Now let $\as(f)$ denote the anti-symmetrization 
\begin{equation}
\label{a_zitra_mne_to_kolo_naprosto_znici}
\as(f)(v_1,\ldots,v_k) 
:= \sum_{\sigma\in \Sigma_k}\chi\cdot f(v_{\sigma(1)},\ldots,v_{\sigma(k)})
\end{equation}
of a linear map $f:V^{\otimes k}\rightarrow V$.

\begin{theorem}
\label{B}
The symmetrization of the (non-symmetric) braces on 
$\bigoplus_{k \geq 1} \Hom(V^{\otimes k},V)$ 
constructed in~\cite{GV} coincides with
the symmetric braces~(\ref{nikdo}). By this we mean that, for each $f,
g_1,\ldots,g_n \in \Hom(V^{\otimes *},V)$,
$$
\sum_{\sigma \in \Sigma_n}
\epsilon \cdot \as(f\{g_{\sigma(1)},\ldots,g_{\sigma(n)}\})=\as(f)\langle
\as(g_1),\ldots,\as(g_n)\rangle,
$$
where $\epsilon$ is the Koszul sign of the permutation
\[
(g_1,\ldots,g_n) \longmapsto (g_{\sigma(1)},\ldots, g_{\sigma(n)}).
\]
\end{theorem}

The proofs of these theorems will appear
in~\cite{daily-progress}. Theorems~\ref{A} and~\ref{B} can be combined
into:

\begin{corollary}
Let us consider the space $\bigoplus_{k \geq 1} \Hom(V^{\otimes
k},V)$ with the symmetric brace algebra structure given by
the symmetrization of the (non-symmetric) brace algebra of~\cite{GV}. Consider also the space of anti-symmetric maps  
$\bigoplus_{k \geq 1} \Hom(V^{\otimes k},V)^{\as}$, with the symmetric
braces~(\ref{nikdo}). Then the anti-symmetrization
$$
\as : \bigoplus_{k \geq 1} \Hom(V^{\otimes k},V)
\to \bigoplus_{k \geq 1} \Hom(V^{\otimes k},V)^{\as}
$$
defined in~(\ref{a_zitra_mne_to_kolo_naprosto_znici}) is a
homomorphism of symmetric brace algebras.
\end{corollary}

As a corollary to Theorem~\ref{B}, we obtain Theorem 3.1 of \cite{LM}:

\begin{corollary}
The anti-symmetrization $l := \as(\mu)$ of an $A_{\infty}$-structure $\mu$
yields an $L_{\infty}$-structure.
\end{corollary}

\begin{proof}  The proof immediately follows from 
$as(\mu\{\mu\})=as(\mu)\langle as(\mu)\rangle=l\langle l\rangle$.
\end{proof}

\begin{remark}
{\rm
There are two degree conventions
for $L_\infty$-algebras. Under the first convention used in~\cite{LM}, 
structure operations are
graded antisymmetric maps $l_k : \otexp Vk \to V$ of degree
$k-2$ which can be pieced together into a degree $-1$ coderivation of
$\coland(\susp V)$. Under the second convention, structure operations are
maps $l_k : \otexp Vk \to V$ of degree $2-k$ which can be encoded by
a degree $+1$ coderivation of $\coland(\desusp V)$. These conventions
differ by the ``flip'' $V_k \mapsto V_{-k}$ of the underlying graded
vector space $V$. In this paper, we use the first convention.
}
\end{remark} 

As we saw in the second half of this section, symmetric brace algebras
may be interpreted as a tool formalizing compositions of
anti-symmetric maps. In the remaining sections, we try to convince the
reader that this formalization can be used for concrete calculations. 

\section{Equivalence between symmetric brace algebras and pre-Lie
         algebras}
\label{sec3}

\def\Tree{{\mathcal T \hskip -.2em {\it ree\/}}} \def\Vert#1{{\it Vert}(#1)}
\def\Edg#1{{\it Edg}(#1)} \def\bfb{{\mathbf b}} \def\bfc{{\mathbf c}}
\def\rrada#1#2#3{#1_{#2},\ldots,#1_{#3}} \def\bfx{{\mathbf x}}
\def\orada#1#2#3{#1_{#2} \otimes \cdots \otimes #1_{#3}}
\def\wphi{{\widetilde \phi}} \def\bfv{{\mathbf v}}
\def\SB{{\rm SB}} \def\pL{{\rm PL}}
\def\opSB{{\mathcal {SB}}}
\def\oppL{{\it p{\mathcal L}\it ie}}
\def\Span{{\it Span}}

In this section we show that the category of symmetric brace algebras
is isomorphic to the category of pre-Lie algebras 
(Proposition~~\ref{pisu_v_Raleigh}). As we already
observed in Remark~\ref{nejak_se_nic_nedari}, this statement was
proved in~\cite{guin-oudom}, but we attempt to give a more direct and
conceptual proof. We show that the free symmetric brace algebra
$\SB(X)$ generated by a set $X$ is functorially isomorphic to the free
pre-Lie algebra $\pL(X)$ generated by the same set. This, by arguments
analyzed for example in~\cite[Section~II.1.4]{markl-shnider-stasheff:book},
means that the operad $\opSB$ describing symmetric brace algebras is
isomorphic to the operad $\oppL$ for pre-Lie algebras, which
tautologically implies that the corresponding {\em categories of
algebras\/} are isomorphic. Since both operads $\opSB$ and $\oppL$ are
concentrated in degree zero, it is enough to consider in this section
only non-graded objects, avoiding thus the sign issue
completely. 

Let us start with some necessary technicalities.  By a {\em rooted
tree\/} we understand a nonempty contractible graph with a
distinguished vertex called the {\em root\/}.  We denote $\Tree$ the
set of all rooted trees. For $T \in \Tree$, let $\Vert T$ (resp.~$\Edg
T$) denote the set of vertices (resp.~edges) of $T$. We will call the
unique rooted tree with one vertex $r$ (which is also its root) the
{\em singleton\/} and denote it $\{r\}$.

Let $S, S_1, \ldots, S_n \in \Tree$, $n \geq 1$, be rooted trees and
let $\bfv : \{\rada S1n\} \to \Vert S$ be a map that assigns to each
element of the set $\{\rada S1n\}$ a vertex of $S$.  Let $S\{\rada
S1n\}_\bfv \in \Tree$ be the tree obtained by connecting, for $1 \leq
i \leq n$, the root of $S_i$ to the vertex $\bfv(S_i) \in \Vert
S$. Therefore
\begin{eqnarray*}
\Vert {S\{\rada S1n\}_\bfv} &=& \Vert S \sqcup \Vert{S_1}
\sqcup \cdots \sqcup \Vert{S_n} \mbox { and }
\\
\Edg {S\{\rada S1n\}_\bfv} &=&  \Edg S \sqcup \Edg{S_1}
\sqcup \cdots \sqcup \Edg{S_n} \sqcup \{\rada e1n\},
\end{eqnarray*}
where $e_i$ is a new edge that joints the vertex $\bfv(S_i) \in \Vert
S$ with the root of the tree~$S_i$. Let us emphasize that in the
notation $S\{S_1,\ldots,S_n\}_\bfv$ the curly braces indicate that the
construction depends only on the {\em set\/} of trees $\rada S1n$ not
on their relative order and have therefore nothing in common with the
non-symmetric braces considered in Section~\ref{sec2}.

\begin{example}
\label{zimni_deprese?}
{\rm\
Let $T$ be a tree with at least two vertices and let $r \in \Vert T$
be its root. Then there exist a unique set $\{T_1,\ldots,T_m\}$ of
rooted trees such that
\[
T = \{r\}\{T_1,\ldots,T_m\}_\bfv,
\]
where $\bfv : \{T_1,\ldots,T_m\} \to \{r\}$ is the unique set map that
sends each $T_i$ to the root $r$ of the singleton $\{r\}$. Since
$\bfv$ carries no information, we will drop it from the notation and
write simply $T = \{r\}\{T_1,\ldots,T_m\}$.
}
\end{example}

Let $B$ be an arbitrary symmetric brace algebra, $T \in \Tree$ a
rooted tree and $\bfb : \Vert T \to B$ a set map. We may interpert
$\bfb$ as a {\em decoration\/} of the vertices of $T$ with elements of
$B$. For such a decoration we define the element $T(\bfb) \in B$
inductively as follows.

If $T = \{r\}$ is the singleton, we put $T(\bfb) := \bfb(r) \in B$.
If $T$ has at least two vertices, it decomposes as in
Example~\ref{zimni_deprese?} into $T = \{r\}\{\rada T1m\}$. Using
the restrictions $\bfb_j := \bfb|_{\Vert {T_j}}: \Vert {T_j}
\to B$, $1 \leq j \leq m$, we define
\begin{equation}
\label{fewer}
T(\bfb) = \{r\}\{\rada T1m\}(\bfb)  
:= \bfb(r) \langle T_1(\bfb_1),\ldots,T_m(\bfb_m)\rangle \in B. 
\end{equation}
In the above display, $-\langle -,\ldots,- \rangle$ denotes the
symmetric brace of $B$.  Since each $T_j$, $1 \leq j \leq m$, has
strictly fewer vertices than $T$, the elements $T_j(\bfb_j) \in B$ have
already been defined by induction.  

In the following proposition, which is the main technical result of
this section, we formulate an extension of
axiom~(\ref{pozitri_prohlidka_a_nejsem_zdrav}) of Definition~\ref{sb}.

\begin{proposition}
\label{kourim_vodni_dymku}
Let $B$ be an arbitrary brace algebra, $S,S_1 , \ldots, S_n \in \Tree$
rooted trees and $\bfc : \Vert S \to B, \bfc_i : \Vert {S_i} \to
B$, $1 \leq i \leq n$, decorations of vertices.  Then
\begin{eqnarray}
\label{ale_prohlidkou_jsem_prosel}
S(\bfc)\langle S_1(\bfc_1),\ldots, S_n(\bfc_n) \rangle =
\sum_{\bfv : \{S_1,\ldots,S_n\} \to \Vert S} 
S\{\rada S1n\}_\bfv(\bfc\sqcup \bfc_1 \sqcup \cdots \sqcup \bfc_n),
\end{eqnarray}
where 
\[
\bfc \sqcup \bfc_1 \sqcup \cdots \sqcup \bfc_n : \Vert
{S\{\rada S1n\}_\bfv} =  \Vert S \sqcup \Vert{S_1}
\sqcup \cdots \sqcup \Vert{S_n} \to B
\] 
is the decoration induced in the obvious
way from the decorations $\bfc,\rada {\bfc}1n$.
\end{proposition}

\noindent 
{\bf Proof.}
The proof is given by induction on the number of vertices of $S$. If
$S$ is the singleton, then~(\ref{ale_prohlidkou_jsem_prosel}) 
follows immediately from the defining equation~(\ref{fewer}).
If $S$ has at least two vertices, then
\[
S = \{r\}\{\rada T1m\}
\]
for some $T_j \in \Tree$ as in Example~\ref{zimni_deprese?}.
Let $\bfb_j := \bfc|_{\Vert {T_j}}$, $1 \leq j \leq m$. Then the left
hand side of~(\ref{ale_prohlidkou_jsem_prosel}) can be expanded into
\[
\bfc(r)\langle T_1(\bfb_1),\ldots,T_m(\bfb_m)\rangle
\langle S_1(\bfc_1),\ldots,S_n(\bfc_n) \rangle.
\]
In the rest of this proof, we simplify the notation by writing
$S_i(\bfc)$ instead of $S_i(\bfc_i)$, $1 \leq i \leq n$, and
$T_j(\bfb)$ instead of $T_j(\bfb_j)$, $1 \leq j \leq m$. We believe
that such a simplification will not confuse the reader. With this
convention assumed, the left hand side
of~(\ref{ale_prohlidkou_jsem_prosel}) reads
\[
\bfc(r)\langle T_1(\bfb),\ldots,T_m(\bfb)\rangle
\langle S_1(\bfc),\ldots,S_n(\bfc) \rangle.
\]
Axiom~(\ref{pozitri_prohlidka_a_nejsem_zdrav}) converts this expression into
\[
\sum_{\it ush}
\bfc(r) 
\left\langle T_1(\bfb)
\langle S_{i^1_1}(\bfc),\squeezedldots, S_{i^1_{t_1}}(\bfc) \rangle,
\ldots,
T_m(\bfb)
\langle S_{i^m_1}(\bfc),\squeezedldots, S_{i^m_{t_m}}(\bfc)
\rangle,
S_{i^{m+1}_1}(\bfc),
\ldots,S_{i^{m+1}_{t_{m+1}}}(\bfc) 
\right\rangle,
\]
where the summation runs over the same set of unshuffles as
in~(\ref{pozitri_prohlidka_a_nejsem_zdrav}). 
Since each $T_j$ has fewer vertices than $S$, we may use the induction
to convert the above expression into
\[
\sum_{\it ush}  
\hskip -.2em   \sum_{\bfv_1,\ldots \bfv_m} \hskip -.2em  \bfc(r)
\left\langle 
T_1\{S_{i^1_1},\squeezedldots, S_{i^1_{t_1}}\}_{\bfv_1}
(\bfb \sqcup \bfc),
\squeezedldots,
T_m\{S_{i^m_1},\squeezedldots, S_{i^m_{t_m}}\}_{\bfv_m}
(\bfb \sqcup \bfc), 
S_{i^{m+1}_1}(\bfc),
\squeezedldots,S_{i^{m+1}_{t_{m+1}}}(\bfc)
\right\rangle
\]
where each $\bfv_j$ runs over all set maps 
$\bfv_j : \{S_{i^j_1},\ldots, S_{i^1_{t_j}}\} \to \Vert {T_j}$, $ 1 \leq j
\leq m$.

Using~(\ref{fewer}), the above display can be written as
\[
\sum_{\it ush} \sum_{\bfv_1,\ldots \bfv_m}
\{r\}\left\{T_1\{S_{i^1_1},\squeezedldots, S_{i^1_{t_1}}\},
\ldots,
T_m\{S_{i^m_1},\squeezedldots, S_{i^m_{t_m}}\},S_{i^{m+1}_1}
\ldots,S_{i^{m+1}_{i_{m+1}}}
\right\}(\bfb \sqcup \bfc)
\]
which can be easily identified with
\[
\sum_{\bfv : \{S_1,\ldots,S_n\} \to \Vert S} 
S\{\rada S1n\}_\bfv(\bfc)
\]
which is the right hand side of~(\ref{ale_prohlidkou_jsem_prosel}).%
\qed

\begin{exercise}
{\em 
Show that axiom~(\ref{pozitri_prohlidka_a_nejsem_zdrav})
is a particular case of~(\ref{ale_prohlidkou_jsem_prosel}) for $S$
the corrola
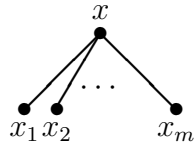
\begin{center}
\unitlength .35em
\linethickness{0.4pt}
\ifx\plotpoint\undefined\newsavebox{\plotpoint}\fi 
\begin{picture}(26,10.5)(0,9)
\thicklines
\put(14,17){\line(-1,-1){7}}
\put(7,10){\line(0,1){0}}
\put(14,17){\line(-3,-5){4.2}}
\put(14,17){\line(1,-1){7}}
\put(21,10){\line(0,1){0}}
\put(7,10){\makebox(0,0)[cc]{$\bullet$}}
\put(21,10){\makebox(0,0)[cc]{$\bullet$}}
\put(10,10){\makebox(0,0)[cc]{$\bullet$}}
\put(14,12){\makebox(0,0)[cc]{$\cdots$}}
\put(14,17){\makebox(0,0)[cc]{$\bullet$}}
\put(7,8){\makebox(0,0)[cc]{$x_1$}}
\put(10,8){\makebox(0,0)[cc]{$x_2$}}
\put(21,8){\makebox(0,0)[cc]{$x_m$}}
\put(14,19){\makebox(0,0)[cc]{$x$}}
\end{picture}
\end{center}
with the root decorated by $x \in B$ and the remaining vertices by
$\rada x1m \in B$, and $S_i$'s the singletons with the unique vertex
decorated by $y_i$, $1 \leq i \leq n$.
}
\end{exercise}

Let us describe our realization of the free symmetric brace algebra. 
For a set $X$ define
\[
\SB (X) : =   \Span\{X_T;\ T \in \Tree\}, 
\]
where $X_T$ is the set of all decorations $\bfx: \Vert T \to X$ of the
vertices of the tree $T$ by elements of $X$. 
Let us introduce symmetric braces on the vector space 
$\SB(X)$ as follows. 
For $\bfx : \Vert S \to X \in X_T$ and  $\bfx_i : \Vert {S_i} \to X
\in X_{S_i}$, $1 \leq i \leq n$, set
\[
\bfx \langle \rada {\bfx}1n \rangle := \sum_{\bfv : \{\rada S1n\}
\to \Vert S} \bfx \sqcup \bfx_1 \sqcup \cdots \sqcup \bfx_n,
\] 
where 
$\bfx \sqcup \bfx_1 \sqcup \cdots \sqcup \bfx_n  
: \Vert {S\{\rada S1n\}_\bfv} \to X$ is the obvious induced decoration
of $S\{\rada S1n\}_\bfv$.
The canonical map $i: X \hookrightarrow \SB (X)$ sends $x \in X$ to the
singleton $\{r\}$ decorated by $x$.  

\begin{theorem}
The object $X \hookrightarrow \SB(X)$ defined above is the free
symmetric brace algebra on~$X$.
\end{theorem}

\noindent 
{\bf Proof.}
Observe first that Propostion~\ref{kourim_vodni_dymku}
implies that  $\SB(X)$ is indeed a symmetric brace algebra.
One must show next that for an arbitrary symmetric brace algebra $B$ and for
an arbitrary set map $\phi : X \to B$ there exist a unique
homomorphism  $\wphi : \SB(X) \to B$ of symmetric brace algebras 
for which the diagram
\begin{center}
\unitlength .6em
\thicklines
\begin{picture}(26,11)(0,6)
\thinlines
\put(6.75,12.75){\bezier{50}(0,0)(0,.25)(.25,.25)\bezier{50}(.25,.25)(.5,.25)(.5,.0)}
\put(9,14){\vector(1,0){11}}
\put(22,14){\makebox(0,0)[cc]{$B$}}
\put(7,14){\makebox(0,0)[cc]{$X$}}
\put(20.4,13.5){\vector(2,1){.07}}
\multiput(8.93,7.93)(.7857,.3571){15}{{$\cdot$}}
\put(7,7){\makebox(0,0)[cc]{$\SB (X)$}}
\put(6.75,12.75){\vector(0,-1){4.75}}
\put(13,14.75){\makebox(0,0)[bc]{$\phi$}}
\put(6,11){\makebox(0,0)[cc]{$i$}}
\put(13,10.75){\makebox(0,0)[bc]{$\widetilde \phi$}}
\end{picture}
\end{center}
commutes.
A moment's reflection convinces us that the only possible choice is
\[
\wphi(\bfx) := T(\phi \circ \bfx),
\] 
where, for $\bfx : \Vert T \to X \in X_T$, $\phi \circ \bfx$ is the
composition $\Vert T \stackrel
{\bfx}{\longrightarrow} X \stackrel {\phi}{\longrightarrow} B$.  Such
a $\wphi : \SB(X) \to B$ is, again by
Proposition~\ref{kourim_vodni_dymku}, indeed a symmetric brace algebra
homomorphism. This finishes the proof.%
\qed

The following proposition is due to D.~Guin and J.-M. Oudom~\cite{guin-oudom}.

\begin{proposition}
\label{pisu_v_Raleigh}
The category of symmetric brace algebras is isomorphic to the category
of pre-Lie algebras.
\end{proposition}

\noindent 
{\bf Proof.}
It is immediate to see that the pre-Lie algebra $\SB(X)_{\it pL}$
associated to the free brace algebra $\SB(X)$ is canonically
isomorphic to the free pre-Lie algebra $\pL(X)$ as described
in~\cite{chapoton-livernet:pre-lie}. Proposition~{II.1.27}
of~\cite{markl-shnider-stasheff:book} implies that the corresponding
operads $\opSB$ and $\oppL$ are isomorphic, and the proposition
follows.%
\qed

\newtheorem{situation}[theorem]{Situation}
\newtheorem{lemma}[theorem]{Lemma}

\def\bfmu{\mbox{$\mu \hskip -2.4mm \mu \hskip - 2.4mm \mu$}}
\def\bfnu{\mbox{$\nu \hskip -2.4mm \nu \hskip - 2.4mm \nu$}}
\def\bfl{\mbox{$l \hskip -1.5mm l \hskip - 1.5mm l$}}
\def\bfk{\mbox{$k \hskip -2.4mm k \hskip - 2.4mm k$}}
\def\p{\mbox{$p \hskip -2.15mm p \hskip - 2.15mm p$}}
\def\pkernel{p-kernel} \def\barmu{{\overline{\mu}}}
\def\pa{\partial} \def\barl{{\overline{l}}} \def\calP{{\mathcal P}}
\def\barnu{{\overline{\nu}}} \def\bark{{\overline{k}}}
 \def\End{{\it End\/}}

\section{Transfers of strongly homotopy structures}
\label{sec4}
In this section we show how brace algebras can be
used to simplify formulas for transfers of strongly homotopy structures.
In~\cite{markl:ha,markl:tr} we considered the following situation.

\begin{situation}
\label{sit}
We are given chain complexes $(V,\pa_V)$, $(W,\pa_W)$ and chain maps
$f : (V,\pa_V) \to (W,\pa_W)$, $g : (W,\pa_W) \to (V,\pa_V)$ such
that the composition $gf$ is chain homotopic to the identity $\id_V:
V\to V$, via a chain-homotopy $h$. In other words, $g$ is a left
homotopy inverse of $f$.
\end{situation}

We assumed that $(V,\pa_V)$ was equipped with an $A_\infty$-structure 
$\bfmu = (\mu_1,\mu_2,\mu_3,\ldots)$ with $\mu_1 = \pa_V$. We were looking for
an $A_\infty$-structure $\bfnu =(\nu_1,\nu_2,\nu_3,\ldots)$ on $(W,\pa_W)$,
with $\nu_1 = \pa_W$, such that the $A_\infty$-structures $\bfmu$ and
$\bfnu$ were equivalent, via a suitable extensions of the chain maps
$f$ and $g$, see~\cite[Problem~2]{markl:tr} for a precise
formulation. 

We call the $A_\infty$-structure $\bfnu$ with the properties specified
in the above paragraph the {\em
transfer\/} of $\bfmu$. The existence of transfers follows from
general principles (see~\cite{markl:ha}), but 
in~\cite{markl:tr} we constructed such 
a  $\bfnu = (\nu_1,\nu_2,\nu_3,\ldots)$ explicitly by the
``Anzatz''
\begin{equation}
\label{An} 
\nu_n : = f \circ \p_n \circ \otexp g n,\ n \geq 2, 
\end{equation}
where $\p_n :
\otexp Vn \to V$ (the {\em \pkernel\/}) was a degree $n-1$ linear
map defined inductively by
$\p_2 : = \mu_2$ and
\begin{equation}
\label{Dny_otevrenych_dveri}
\p_n := \sum_{B} (-1)^{\vartheta(r_1,\ldots,r_k)}
\mu_k ( h \circ \p_{r_1} \otimes \cdots \otimes h \circ \p_{r_k}),
\ \mbox { for $n \geq 2$}. 
\end{equation}
In the above display we denoted
\[
\vartheta(u_1,\ldots,u_s) := \sum_{1 \leq \alpha < \beta \leq s}
u_\alpha(u_\beta + 1),
\]
the summation was taken over
\[
B  :=  \{k, \Rada r1k\ |\ 2 \leq k \leq n,\ 
r_1,\ldots, r_k \geq 1,\ r_1 + \cdots + r_k = n\}
\]
and the formal convention that $h \p_1 = \id$ was assumed.
We proved that the \pkernel\ satisfies
\begin{equation}
\label{Phily}
\pa(\p_n)
=
\sum_A
(-1)^{i(l+1) +n} \p_{k}(\id^{\ot i-1} \ot  gf \circ \p_l \ot
\id^{\ot k-i}),\ n \geq 2,
\end{equation}
where $\pa$ is the differential on $\Hom(\otexp Vn,V)$ induced from
$\pa_V$ in the standard way. We then
derived from~(\ref{Phily}) that the operations defined by the
Anzatz~(\ref{An}) indeed form an $A_\infty$-structure on $(W,\pa_W)$.

Let us translate the above calculations into the language of brace
algebras. 
Consider the graded vector space $\End_*(V)$ defined by
\[
\End_s(V) := 
\bigoplus_{p-k+1=s}\Hom(V^{\otimes k},V)_p
\]
with the (nonsymmetric) brace algebra structure introduced in~\cite{GV}. 
As we observed in
Example~\ref{Jsem_zvedav_jestli_se_ozve_pred_Sevillou},  
the $A_\infty$-structure $\bfmu =
(\mu_1,\mu_2, \mu_3,\ldots)$ assembles into an element $\mu := \mu_1 + \mu_2 +
\mu_3 + \cdots \in
\End_{-1}(V)$ that satisfies $\mu\{\mu\} = 0$. Let $\barmu := \mu_2 +
\mu_3 + \cdots$ so that $\mu = \pa_V + \barmu$ (recall that we
assumed $\mu_1 = \pa_V$). Therefore
\begin{equation}
\label{zitra_jedu_do_Lancasteru}
\mu\{\mu\} = \pa_V\{\pa_V\} + \pa_V\{\barmu\} + \barmu\{\pa_V\} + 
\barmu\{\barmu\} = 0
\end{equation}
where, of course, $\pa_V\{\pa_V\} =0$.
Let, for $u \in \End_*(V)$, 
\[
\pa u := \pa_V\{u\} - (-1)^{|u|} u
\{\pa_V\}.
\] 
With this notation assumed, 
the defining identity~(\ref{zitra_jedu_do_Lancasteru}) for 
$A_\infty$-structures can be rewritten as
$\pa \barmu + \barmu\{\barmu\} = 0.$
The same analysis takes place also for $\bfnu =
(\nu_1,\nu_2,\nu_3,\ldots)$, 
that is, $\barnu :=
\nu_2 + \nu_3 + \cdots \in \End_{-1}(W)$ defines an
$A_\infty$-structure on $(W,\pa_W)$ if and only if 
\begin{equation}
\label{master}
\pa \barnu + 
\barnu\{\barnu\} = 0.  
\end{equation}

The maps $f: V \to W$ and $g: W \to V$ induce
a degree zero map $\Phi : \End_*(V) \to \End_*(W)$ that
sends $u \in \Hom(\otexp Vn,V)$ to $f \circ u \circ \otexp gn$. 
It is clear that $\pa \Phi = \Phi \pa$ and that
\begin{equation}
\label{koupil_jsem_Messiana}
\Phi(u)\{\Phi(u_1),\ldots,\Phi(u_n)\} =
\Phi(u\{gf \circ u_1,\ldots,gf \circ u_n\}),
\end{equation}
for $u, \rada u1n \in \End_*(V)$.

Let finally $\p := \p_2 + \p_3 + \cdots \in \End_{-1}(V)$. In
the above notation, the Anzatz~(\ref{An}) reads:
\begin{equation}
\label{An1}
\barnu := \Phi(\p).
\end{equation}
With a little effort, one verifies that the 
inductive formula~(\ref{Dny_otevrenych_dveri}) defining the \pkernel\
can be rewritten as
\begin{equation}
\label{zitra_jedu_za_MF}
\p = \barmu + \barmu\{h\circ \p\} + \barmu \{h\circ\p,h\circ\p\} + \barmu
\{h\circ\p,h\circ\p,h\circ\p\} + \cdots.
\end{equation}
Setting formally $\mu\{1\} := \mu$, the above display 
can further be shortened into
\[
\p = \barmu\left\{\frac{1}{1 - h \circ\p} \right\}.
\]
In the same vein, formula~(\ref{Phily}) reads:
\begin{equation}
\label{Phily1}
\pa \p + \p\{gf\circ \p\} = 0.
\end{equation}

Let us show that~(\ref{Phily1}) indeed implies the master
equation~(\ref{master}) for $\barnu$ defined 
by~(\ref{An1}).
We have
\[
\pa \barnu = \pa \Phi(\pa \p) = \Phi(\pa \p) = - \Phi(\p \{gf \circ
\p\}) = - \Phi(\p)\{\Phi (\p)\} = - \barnu\{\barnu\}.
\]
where we used~(\ref{koupil_jsem_Messiana}) and the fact that $\Phi$
is a chain map.
Therefore $\pa \barnu =- \barnu\{\barnu\}$, which is~(\ref{master}). 
We leave as an exercise to
prove~(\ref{Phily1}) using axioms of brace algebras only.

Let us consider the $L_\infty$-version of the above problem.  
In Situation~\ref{sit}, we are given an $L_\infty$-structure 
$\bfl = (l_1,l_2,l_3,\ldots)$, with $l_1 = \pa_V$, on $(V,\pa_V)$. 
We are looking for
its transfer $\bfk =(k_1,k_2,k_3,\ldots)$ onto $(W,\pa_W)$, with $k_1 = \pa_W$.

It is clear that instead of $\End_*(V)$ which helped us with the
$A_\infty$-case, 
we need the symmetric brace algebra
$B_*(V)$ introduced in 
Example~\ref{Jsem_zvedav_jestli_se_ozve_pred_Sevillou}. Let $\barl := l_2 + l_3 +
\cdots \in B_{-1}(V)$. The following theorem gives an explicit
construction of a transfer $\bfk$ of $\bfl$.

\begin{theorem}
\label{ii}
The \pkernel\ $\p = \p_2 + \p_3 + \cdots \in B_{-1}(V)$ 
defined inductively by $\p_2
:= l_2$ and
\begin{equation}
\label{jakpak_je_doma?}
\p := \barl + \barl\angles {h \circ \p} + \frac1{2!} 
\barl \angles {h\circ \p,h\circ \p} +
\frac1{3!} \barl \angles {h\circ \p,h\circ \p,h\circ \p} + \cdots
\end{equation}
satisfies
\begin{equation}
\label{je_doma_vse_v_poradku?}
\pa \p + \p \angles { gf \circ\p} =0.
\end{equation}
The Anzatz
\begin{equation}
\label{An2}
\bark := \Phi(\p)
\end{equation}
then defines an $L_\infty$-structure $\bfk = (k_1,k_2,k_3,\ldots)$,
with $k_1 = \pa_W$ on $(W,\pa_W)$, where $\Phi$ is defined in exactly the same way as in the $A_\infty$-case.
\end{theorem}

Equation~(\ref{jakpak_je_doma?}) is the symmetrization in the sense of Theorem \ref{B}
of~(\ref{zitra_jedu_za_MF}). We will shorten it into
\begin{equation}
\label{a_co_zviratka?}
\p = \barl \angles {\exp (h\circ\p)}.
\end{equation}
Our proof of Theorem~\ref{ii} will be based on the following lemma.

\begin{lemma}
\label{byl_jsem_v_Korejske_restauraci}
Elements $a,b,c$,  with
$\deg(b) = \deg(c) = -1$, of an arbitrary symmetric brace algebra satisfy 
\begin{equation}
\label{aa}
a \angles {\exp c} \angles b = a \angles {c \angles b,\exp c} + a
\angles {b,\exp c}
\end{equation}
and
\begin{equation}
\label{bb}
a \angles { b} \angles {\exp c} = a \angles { b\angles{\exp c},\exp c}.
\end{equation}
\end{lemma}

\noindent
{\bf Proof.}
The left hand side of~(\ref{aa}) reads
\begin{eqnarray*}
a \angles {\exp c} \angles b &=&
a \angles {b} + a \angles {c}\angles {b} 
+ \frac1{2!} a \angles {c,c}\angles {b} 
+ \frac1{3!} a \angles {c,c,c}\angles {b} + \cdots
\\
&=& a \angles {b}
\\ 
&& + \hskip .2em a  \angles {b,c} + a\angles {c\angles {b}}
\\ 
&& + \hskip .2em\frac 1{2!}\left( a \angles {b,c,c} + 2a \angles {c \angles
  {b},c} \right)
\\
&& +\hskip .2em \frac 1{3!}\left( a \angles {b,c,c,c} + 3a \angles {c \angles
  {b},c,c} \right) + \cdots
\\
&=&
a \angles {b} + a \angles {b,c} +  \frac 1{2!} a \angles {b,c,c} +  
\frac 1{3!} a \angles {b,c,c,c} + \cdots
\\
&& + \hskip .2em a\angles {c\angles {b}} + a \angles {c \angles {b},c} + 
\frac1{2!} a \angles {c \angles {b},c,c} +
\frac1{3!} a \angles {c \angles {b},c,c,c} + \cdots
\\
&=&  a \angles {c \angles b,\exp c} + a
\angles {b,\exp c}.
\end{eqnarray*}
This proves~(\ref{aa}). Equation~(\ref{bb}) can be proved similarly.
\qed

\noindent
{\bf Proof of Theorem~\ref{ii}.}
Let us apply the differential $\pa$ to the \pkernel\ $\p$ 
defined by~(\ref{a_co_zviratka?}).
We obtain
\begin{eqnarray}
\label{koupil_jsem_si_ryzovy_puding}
\lefteqn{
\pa \p = - \barl \angles \barl \angles {\exp (h\circ \p)}
+ \barl \angles {\p,\exp (h\circ \p)} } \hskip 2em
\\ \nonumber 
&& \hskip 2em
- \barl \angles {gf\circ  \p,\exp (h\circ\p)} 
- \barl \angles {h\circ \p \angles {gf\circ \p},\exp (h\circ\p)} 
\end{eqnarray}
In the above equation we used $\pa \bar l = \barl \angles {\barl}$ and
the high school formula
\[
\pa(\exp (h \circ \p)) = \pa(h \circ \p) \exp (h \circ \p),
\] 
where, of course,
\[
\pa(h \circ \p) = \pa h \circ \p + h \circ \pa \p,
\]
with $\pa h = \id - gf$ and $\pa \p = - \p \angles { gf \circ \p}$ by induction.

By the definition~(\ref{a_co_zviratka?}) of $\p$, the
second term in the right hand side of~(\ref{koupil_jsem_si_ryzovy_puding}) equals 
\[
\barl \angles {\barl \angles {\exp (h\circ\p)}, \exp (h\circ\p) },
\]
while equation~(\ref{bb}) of
Lemma~\ref{byl_jsem_v_Korejske_restauraci}, 
with $a = b = \barl$ and $c = h \circ \p$, gives
\[
-
\barl \angles {\bar l} \angles {\exp (h \circ p)} 
+\barl \angles {\barl \angles {\exp (h\circ\p)}, \exp (h\circ\p) } =0.
\]
Therefore 
the first two terms in the right hand side 
of~(\ref{koupil_jsem_si_ryzovy_puding}) cancel. 
The remaining two terms 
combine, by equation~(\ref{aa}) of
Lemma~\ref{byl_jsem_v_Korejske_restauraci} 
with $a = \barl$, $b = gf \circ \p$ and $c = h \circ \p$, 
into $-\barl\angles { \exp (h \circ p)} \angles { gf \circ \p}$.
Since $\barl \angles { \exp (h \circ \p)} = \p$ by~(\ref{a_co_zviratka?}), 
equation~(\ref{koupil_jsem_si_ryzovy_puding}) implies that 
$\pa \p = - \p \angles { gf \circ \p}$, which
is~(\ref{je_doma_vse_v_poradku?}). 
To prove that $\bark$ defined by the Anzatz~(\ref{An2})
fulfills $\bark \angles {\bark} = 0$ is equally simple
as in the $A_\infty$-case discussed in the first part of this section.%
\qed

\hskip .3em

Although we know from Proposition~\ref{pisu_v_Raleigh} 
that any symmetric brace algebra is
generated by its pre-Lie part,  it is unclear how to
write the defining formula~(\ref{jakpak_je_doma?}) 
for $\p$ using the pre-Lie multiplication only. 
We believe that this demonstrates that symmetric brace algebras are
a useful concept even if they are ``formally'' the same as
pre-Lie algebras.
Since the operadic cochain complex of an arbitrary algebra
over a quadratic Koszul operad $\calP$ carries a natural symmetric brace
algebra structure, the formulas of Theorem~\ref{ii} in fact
define transfers for {\em arbitrary\/} strongly homotopy $\calP$-algebras.

\section{The $L_{\infty}$-structure of a Gauge Algebra}
\label{sec5}

In this section we provide an example of how symmetric brace algebras may be
 used to organize the consequences of the algebraic assumptions that lead to 
 a particular $L_{\infty}$-algebra structure.  This algebraic data, originally described by Berends, Burgers, and van Dam in their analysis of a particular type of gauge theory \cite{BBvD}, was
 recast in \cite{FLS} and was shown to lead to an $L_{\infty}$-algebra structure on the
 space of gauge parameters together with the space of fields.

We summarize this situation in the following fashion. 
Let $\Xi$ denote the (nongraded) vector space of
gauge parameters and $\Phi$ denote the (nongraded) vector space of fields. 
The ``action" is given by a gauge
transformation which is phrased as a linear map 
$\delta:\Xi\rightarrow \Hom(S^*(\Phi),\Phi)$
where $S^*(\Phi)$ is the cofree nilpotent cocommutative coassociative
coalgebra cogenerated by $\Phi$, which in fact coincides with the
linear dual of the algebra of polynomial functions on $\Phi$,
see~\cite[Example~II.3.79]{markl-shnider-stasheff:book}. 
Let us remark that $S^*(\Phi)$ was,
in~\cite{FLS}, denoted $\coland^* \Phi$. This notation was
formally correct, but we think it could be easily mistaken for the
exterior (Grassmann) coalgebra cogenerated by $\Phi$. 
The map $\delta$ is then
extended to a map 
$\hat\delta:\Hom(S^*(\Phi),\Xi)\rightarrow \Hom(S^*(\Phi),\Phi)$ via
\[
\hat\delta(\pi): =ev\circ [(\delta\circ \pi)\otimes \id]\circ\Delta
\] 
where $ev: \Hom (S^* (\Phi),\Phi) \ot S^*(\Phi) \to \Phi$ 
is the evaluation map and
$\Delta$ is the comultiplication on $S^*(\Phi)$.  
Recall that the vector space
$\Hom(S^*(\Phi),\Phi)$ has a canonical Lie 
bracket given by $[f,g]=f\circ \bar{g}-g\circ\bar{f}$
where $\bar{f}$ denotes the extension of a 
linear map $f\in \Hom(S^*(\Phi),\Phi)$ to a
coderivation $\bar{f}$ on $S^*(\Phi)$.

Another ingredient in this picture is the 
assumed existence of a map $C:\Xi\otimes\Xi\rightarrow
\Hom(S^*(\Phi),\Xi)$ that satisfies the (BBvD) hypothesis
$$[\delta(\xi),\delta(\eta)]=\hat\delta C(\xi,\eta)\in 
\Hom(S^*(\Phi),\Phi)$$
for all $\xi,\eta\in\Xi$.  After this map is extended to a map
$\hat{C}:\Hom(S^*(\Phi),\Xi)\otimes 
\Hom(S^*(\Phi),\Xi)\rightarrow \Hom(S^*(\Phi),\Xi)$,
\begin{equation}
\label{prece_jenom_napsala}
\hat{C}(\pi_1,\pi_2)=
ev \circ [(C \ot \id) \circ (\pi_1 \ot \pi_2 \ot \id)] \circ (\Delta \ot
\id) \circ \Delta,
\end{equation}
a bracket that satisfies the Jacobi identity may
be imposed on the space $\Hom(S^*(\Phi),\Xi)$ via
$$[\pi_1,\pi_2]=\pi_1\circ\overline{\hat\delta(\pi_2)}-\pi_2
\circ\overline{\hat\delta(\pi_1)}+
\hat{C}(\pi_1,\pi_2).$$

Next, we consider the graded vector space $\mathbb{L}$ 
with $\mathbb{L}_0:=\Xi$, $\mathbb{L}_{-1}:=\Phi$,
and $\mathbb{L}_n:=0$ for $n\neq0,-1$.  
By Theorem 2 of \cite{FLS}, an $L_{\infty}$-structure may be
defined on $\mathbb{L}$ by constructing a degree $-1$ 
map $D:\coland^*(\susp\mathbb{L})\rightarrow\uparrow
\mathbb{L}$ by piecing together the maps $\delta$ and $C$.  
The Jacobi identity for the bracket on $\Hom(S^*(\Phi),\Xi)$ 
implies that $D\circ\overline{D}=0$, the
$L_{\infty}$ relations for $\mathbb{L}$.

In the remainder of this note, 
we show how this $L_{\infty}$-structure may be explained in terms 
of the symmetric brace algebra $B_*(\mathbb{L})$ on 
$\Hom(\mathbb{L}^{\otimes *},\mathbb{L})^{as}$ constructed in
Example~\ref{Jsem_zvedav_jestli_se_ozve_pred_Sevillou}.   
Specifically, we will use the
symmetric brace algebra structure to construct a 
bracket on $\Hom(\mathbb{L}^{\otimes *},\mathbb{L})^{as}$ 
so that the existence of an
$L_{\infty}$-structure on $\mathbb{L}$ 
will be equivalent to this bracket satisfying the Jacobi identity
on the subspace $\Hom(\Phi^{\otimes *},\Xi)^{as}$.

Let us fix two maps, 
$\nabla$ and $\Upsilon$ in 
$\Hom(\mathbb{L}^{\otimes *},\mathbb{L})^{as}$.  We require that 
\begin{itemize}
\item[(i)]
the map $\nabla$  
have values in $\Phi$ and be the zero map
when the number of inputs from the space $\Xi$ is not equal to 1.  
Similarly, we require that
\item[(ii)]
the map $\Upsilon$ take values in $\Xi$ and be
the zero map when the number of inputs from $\Xi$ is not equal to 2.
\end{itemize}

An important example of maps with the above properties is provided by the
maps $\delta$ and $C$, via the exponential correspondence
\[
\Hom(\Xi, \Hom (S^*(\Phi),\Phi)) \ni \delta \longleftrightarrow 
\nabla \in \Hom(\Xi \ot S^*(\Phi),\Phi)
\]
and
\[
\Hom(\Xi \wedge \Xi, \Hom (S^*(\Phi),\Xi)) \ni C \longleftrightarrow 
\Upsilon \in \Hom(\Xi\wedge \Xi \ot S^*(\Phi),\Xi),
\]
where $\wedge$ denotes the antisymmetric (exterior) product.

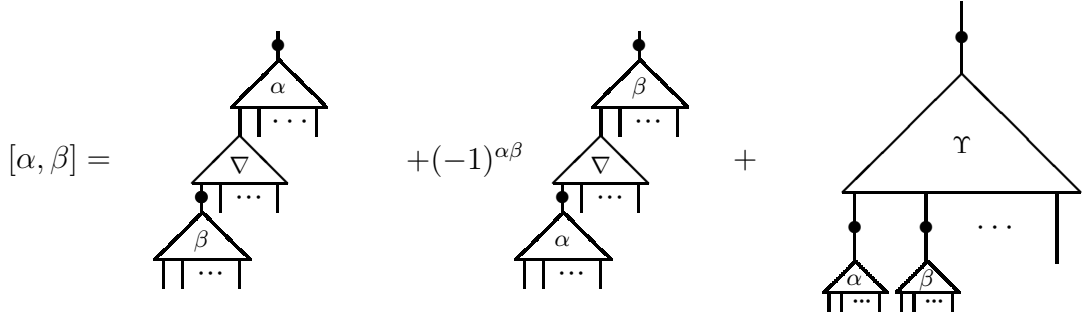
\begin{figure}[ht]
\begin{center}
\unitlength=1.2cm
\begin{picture}(4,4)(3.5,-1.6)
\put(0,0){\makebox(0,0){$[\alpha,\beta] =$}}
\put(2.425,.845){\makebox(0,0){\Alpha{.36}}}
\put(2,0){\makebox(0,0){\NABLA{.36}}}
\put(1.575,-.845){\makebox(0,0){\Beta{.36}}}
\put(4.5,0){\makebox(0,0){$+(-1)^{\alpha \beta}$}}
\put(4,0){
\put(2.425,.845){\makebox(0,0){\Beta{.36}}}
\put(2,0){\makebox(0,0){\NABLA{.36}}}
\put(1.575,-.845){\makebox(0,0){\Alphanobul{.36}}}
}
\put(7.6,0){\makebox(0,0){$+$}}
\put(0,-.7){
\put(10,1){\makebox(0,0){\UPSILON{.9}}}
\put(8.81,-.6){\makebox(0,0){\Alphanobul{.24}}}
\put(9.61,-.6){\makebox(0,0){\Betanobul{.24}}}
}
\end{picture}
\end{center}
\caption{\label{fig:1}
A pictorial description of $[\alpha,\beta]$. Multilinear maps are
symbolized by triangles, with the bullet 
$\bullet$ marking outputs/inputs in $\Xi$.}
\end{figure}

It follows immediately from (i) and (ii) above that $\nabla,\Upsilon
\in B_{-1}(\mathbb{L})$. 
For $\alpha,\beta\in B_*(\bbL)$, we define a degree $-1$ bracket
$$
[\alpha,\beta]=\alpha\langle\nabla\langle\beta\rangle\rangle
+(-1)^{\alpha\beta}\beta\langle\nabla
\langle\alpha\rangle\rangle+\Upsilon\langle\alpha,\beta\rangle.
$$
A pictorial definition of this bracket is given in Figure~\ref{fig:1}. 
We have

\begin{theorem}
The above bracket, restricted to 
$\Hom(\Phi^{\otimes *}, \Xi)^{\rm as} \subset B_1(V)$, satisfies the Jacobi
identity if and only if
\begin{equation}
\nabla\langle\nabla\rangle+\nabla\langle\Upsilon\rangle=0  \label{bbvdhyp}
\end{equation}
and
\begin{equation}
\Upsilon\langle\nabla\rangle+\Upsilon\langle\Upsilon\rangle=0.  \label{jac}
\end{equation}
\end{theorem}
\begin{proof}
When we iterate the bracket, we obtain the expression
\begin{eqnarray*}
[\alpha,[\beta,\gamma]]
&=&\alpha\langle\nabla\langle\beta\langle\nabla\langle\gamma\rangle\rangle
\rangle\rangle-\alpha\langle\nabla\langle\gamma\langle\nabla\langle
\beta\rangle\rangle\rangle\rangle
+\alpha\langle\nabla\langle\Upsilon\langle\beta,\gamma\rangle\rangle\rangle
\\
&&-\beta\langle\nabla\langle\gamma\rangle\rangle
\langle\nabla\langle\alpha\rangle\rangle
+\gamma\langle\nabla\langle\beta\rangle\rangle\langle\nabla
\langle\alpha\rangle\rangle
-\Upsilon\langle\beta,\gamma\rangle\langle\nabla\langle\alpha\rangle\rangle
\\
&&+\Upsilon\langle\alpha,\beta\langle\nabla\langle\gamma\rangle\rangle\rangle
-\Upsilon\langle\alpha,\gamma\langle\nabla\langle\beta\rangle\rangle\rangle
+\Upsilon\langle\alpha,\Upsilon\langle\beta,\gamma\rangle\rangle
\end{eqnarray*}
plus the corresponding terms with $\alpha,\beta,\gamma$ cyclicly
permuted.  
There are two situations
to consider.

We first examine the terms that have $\alpha$ outside of the braces. 
Two such terms that contain two 
$\nabla$'s are
exhibited explicitly above while two more may be found in the cyclic 
permutations of the above expression.  
We let
$$
A=\alpha\langle\nabla\langle\beta\langle\nabla\langle\gamma\rangle\rangle
\rangle\rangle-\alpha\langle\nabla\langle\gamma\langle
\nabla\langle\beta\rangle\rangle\rangle\rangle
-\alpha\langle\nabla\langle\beta\rangle\rangle\langle
\nabla\langle\gamma\rangle\rangle
+\alpha\langle\nabla\langle\gamma\rangle\rangle\langle
\nabla\langle\beta\rangle\rangle
$$

We apply the second relation in the definition of 
symmetric brace algebra twice to the third and fourth terms of $A$
to obtain
$$-\alpha\langle\nabla\langle\beta\rangle
\rangle\langle\nabla\langle\gamma\rangle\rangle=
-\alpha\langle\nabla\langle\beta\rangle,\nabla\langle\gamma\rangle\rangle
-\alpha\langle\nabla\langle\beta,\nabla\langle\gamma\rangle\rangle\rangle
-\alpha\langle\nabla\langle\beta\langle
\nabla\langle\gamma\rangle\rangle\rangle\rangle$$
and
$$\alpha\langle\nabla\langle\gamma
\rangle\rangle\langle\nabla\langle\beta\rangle\rangle=
\alpha\langle\nabla\langle\gamma\rangle,\nabla\langle\beta\rangle\rangle
+\alpha\langle\nabla\langle\gamma,\nabla\langle\beta\rangle\rangle\rangle
+\alpha\langle\nabla\langle\gamma\langle
\nabla\langle\beta\rangle\rangle\rangle\rangle.$$
After adding these to the first two terms of $A$, we have 
$$
A=\alpha\langle\nabla\langle\nabla\langle\beta\rangle,\gamma\rangle\rangle-
\alpha\langle\nabla\langle\nabla\langle\gamma\rangle,\beta\rangle\rangle.
$$
Consider the symmetric algebra relation
\[
\nabla \angles {\nabla} \angles{\beta,\gamma} = 
\nabla \angles {\nabla \angles {\beta,\gamma}} + 
\nabla \angles {\nabla \angles {\beta} ,\gamma} -
\nabla \angles {\nabla \angles {\gamma},\beta} +
\nabla \angles{\nabla,\beta,\gamma}
\]
in which the first and the last summands are equal 
to zero because in each, the map $\nabla$ has more than one
input from $\Xi$.  Consequently,
\begin{equation}
\label{napise_mi_v_dohledne_dobe?}
\nabla \angles {\nabla} \angles{\beta,\gamma} = 
\nabla \angles {\nabla \angles {\beta} ,\gamma} -
\nabla \angles {\nabla \angles {\gamma},\beta},
\end{equation}
therefore
\[
A =\alpha\langle\nabla\langle\nabla\rangle\langle\beta,\gamma\rangle\rangle.
\]
Equation~(\ref{napise_mi_v_dohledne_dobe?}) is illustrated in
Figure~\ref{fig:2}.
\begin{figure}[t]
\begin{center}
\unitlength=1.50cm
\begin{picture}(4,4.4)(2,-1.2)
\put(0,-.4){
\put(0,2){\makebox(0,0){\NABLAmod{.9}}}
\put(-.95,.5){\makebox(0,0){\NABLA{.3}}}
}
\put(-.1,-.5){\makebox(0,0){\line(1,0){2.5}}}
\put(-.8,-1){\makebox(0,0){\Beta{.3}}}
\put(.3,-1){\makebox(0,0){\Gamma{.3}}}
\put(1.8,0){\makebox(0,0){$=$}}
\put(-.4,-1){
\put(4.3,2){\makebox(0,0){\NABLAmod{.9}}}
\put(3.35,.7){\makebox(0,0){\NABLAnobul{.24}}}
\put(3.13,.1){\makebox(0,0){\Beta{.24}}}
\put(3.98,.7){\makebox(0,0){\Gammanobul{.24}}}
}
\put(5.9,0){\makebox(0,0){$-$}}
\put(3.7,-1){
\put(4.3,2){\makebox(0,0){\NABLAmod{.9}}}
\put(3.35,.7){\makebox(0,0){\NABLAnobul{.24}}}
\put(3.13,.1){\makebox(0,0){\Gamma{.24}}}
\put(3.98,.7){\makebox(0,0){\Betanobul{.24}}}
}

\end{picture}
\end{center}
\caption{\label{fig:2}
Equation $\nabla \angles {\nabla} \angles{\beta,\gamma} = 
\nabla \angles {\nabla \angles {\beta} ,\gamma} -
\nabla \angles {\nabla \angles {\gamma},\beta}$.
The horizontal line in the left hand side indicates that $\nabla
\angles {\nabla}$ is performed first.
}
\end{figure}
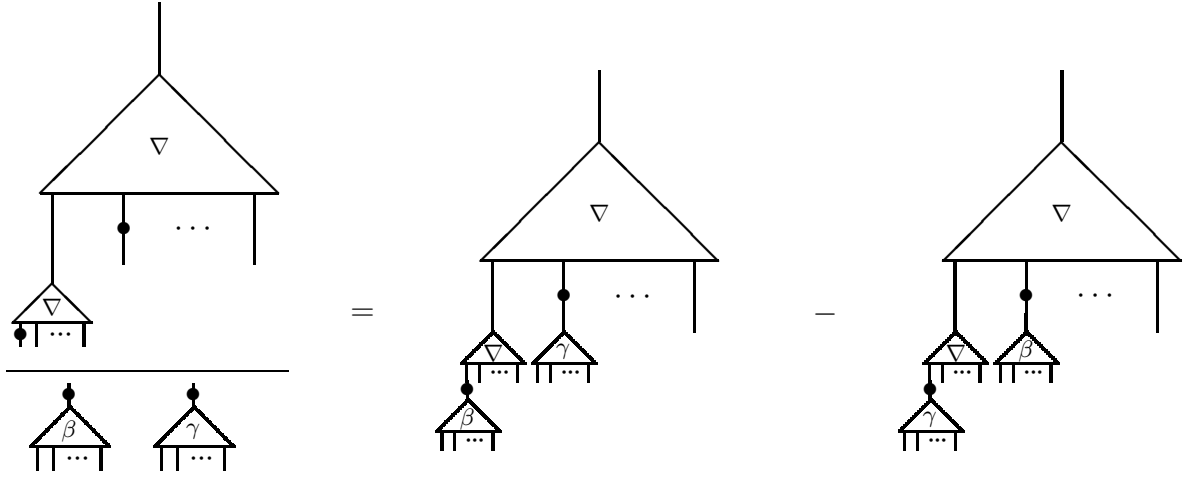  
The remaining term with $\alpha$ outside of 
the braces, $\alpha\langle\nabla\langle\Upsilon\langle\beta,\gamma
\rangle\rangle\rangle$, may be replaced by 
$\alpha\langle\nabla\langle\Upsilon\rangle\langle\beta,\gamma\rangle\rangle$
because by using the symmetric brace relation, we have
$$\alpha\langle\nabla\langle\Upsilon\rangle\langle\beta,\gamma\rangle\rangle
=\alpha\langle\nabla\langle\Upsilon\langle\beta,\gamma
\rangle\rangle\rangle
+\alpha\langle\nabla\langle\Upsilon\langle\beta\rangle,\gamma\rangle\rangle
-\alpha\langle\nabla\langle\Upsilon\langle\gamma\rangle,\beta\rangle\rangle
+\alpha\langle\nabla\langle\Upsilon,\beta,\gamma\rangle\rangle$$
in which the last three terms are equal to zero, 
because in each, the map $\nabla$ has more than one input
from $\Xi$.  When we add this term to $A$, we have
\begin{equation}\alpha\langle\nabla
\langle\nabla\rangle\langle\beta,\gamma\rangle\rangle+
\alpha\langle\nabla\langle\Upsilon\rangle\langle\beta,
\gamma\rangle\rangle  \label{n2}
\end{equation}
together with the similar terms arising from the cyclic 
permutations of $\alpha,\beta$, and $\gamma$.

To account for the remaining terms, we use arguments 
similar to those above to rewrite them as
\begin{equation}
-\Upsilon\langle\nabla\rangle\langle\alpha,\beta,\gamma\rangle
-\Upsilon\langle\Upsilon\rangle\langle\alpha,\beta,\gamma\rangle. \label{nunu}
\end{equation}

Finally, we see that the Jacobi expression 
is equal to zero
if both (\ref{n2}) and (\ref{nunu}) are equal to zero, for all
$\alpha,\beta$ and $\gamma$. It is not
difficult to see that this happens if and only if~(\ref{bbvdhyp})
and~(\ref{jac}) are satisfied (there are ``enough test functions'').
\end{proof}

We note that equation (\ref{bbvdhyp}) now plays the role of the 
(BBvD) hypothesis mentioned above. The following corollary is
basically Theorem~2 of~\cite{FLS}.

\begin{corollary}
If the bracket~(\ref{prece_jenom_napsala}) satisfies the Jacobi identity, 
then the map $l=\nabla+\Upsilon$ is an $L_{\infty}$-%
structure for $\mathbb{L}$.
\end{corollary}

\begin{proof}
As we know from Exercise~\ref{Lowecraft}, we must prove that $l\langle
l \rangle = 0$. But this easily follows from
\[
l \langle l \rangle = (\nabla + \Upsilon)\langle\nabla + \Upsilon
\rangle
= (\nabla\langle \nabla \rangle + \nabla \langle \Upsilon \rangle) +
(\Upsilon\langle \nabla \rangle + \Upsilon \langle \Upsilon \rangle)
\]
and equations~(\ref{bbvdhyp}) and~(\ref{jac}).
\end{proof}


\end{document}